\documentclass[10pt]{amsart}  
\usepackage{amsfonts}
\usepackage{graphicx}
\usepackage{amscd}
\usepackage{amsmath,amsfonts,amssymb,amsthm,mathrsfs,xypic}
\xyoption{curve}
\usepackage{fancybox}
\newcommand{\m}{\Lambda}
\newcommand{\Hom}{\operatorname{Hom}}
\newcommand{\Ker}{\operatorname{Ker}}
\newcommand{\cok}{\operatorname{Coker}}
\newcommand{\Ima}{\operatorname{Im}}
\newcommand{\To}{\operatorname{\longrightarrow}}
\newcommand{\ha}{\operatorname{\mathcal{A}}}
\newcommand {\hp}{\mathcal P}
\newcommand{\s}{\hfill \blacksquare}
\newtheorem{thm}{Theorem}[section]

\newtheorem{cor}[thm]{Corollary}

\newtheorem{lem}[thm]{Lemma}
\newtheorem{exm}[thm]{Example}

\newtheorem{prop}[thm]{Proposition}

\newtheorem{defn}[thm]{Definition}
\begin{document}

\title [Monic representations]
{Monic representations \\ and Gorenstein-projective modules}
\author [Xiu-Hua Luo, Pu Zhang ] {Xiu-Hua Luo, Pu Zhang$^*$}
\thanks{$^*$The corresponding author.}
\thanks{Supported by the NSF of China (10725104).}
\thanks{X.H.Luo$\symbol{64}$sjtu.edu.cn  \ \ \ \ pzhang$\symbol{64}$sjtu.edu.cn}
\dedicatory {Dedicated to the memory of Hua Feng} \maketitle

\begin{center}
Department of Mathematics, \ \ Shanghai Jiao Tong University\\
Shanghai 200240, P. R. China
\end{center}

\begin{abstract} \ Let $\Lambda$ be the path algebra of a finite quiver $Q$
over a finite-dimensional algebra $A$. Then $\Lambda$-modules are
identified with representations of $Q$ over $A$. This yields the
notion of monic representations of $Q$ over $A$. If $Q$ is acyclic,
then the Gorenstein-projective $\m$-modules can be explicitly
determined via the monic representations. As an application, $A$ is
self-injective if and only if the Gorenstein-projective $\m$-modules
are exactly the monic representations of $Q$ over $A$.

\vskip5pt

{\it Key words and phrases. \ representations of a quiver over an
algebra, \ monic representations, \ Gorenstein-projective
modules}\end{abstract}

\section {\bf Introduction}

\vskip10pt

Let $A$ be an Artin algebra, and $A$-mod the category of finitely
generated left $A$-modules. {\it A complete $A$-projective
resolution} is an exact sequence of finitely generated projective
$A$-modules
$$P^\bullet = \cdots \longrightarrow P^{-1}\longrightarrow P^{0}
\stackrel{d^0}{\longrightarrow} P^{1}\longrightarrow \cdots$$ such
that ${\rm Hom}_A(P^\bullet, A)$ remains to be exact. A module $M\in
A$-mod is {\it Gorenstein-projective}, if there exists a complete
$A$-projective resolution $P^\bullet$ such that $M\cong
\operatorname{Ker}d^0$. Let $\mathcal {P}(A)$ be the full
subcategory of $A$-mod of projective modules, and $\mathcal {GP}(A)$
the full subcategory of $A$-mod
 of Gorenstein-projective modules. Then $\mathcal {P}(A) \subseteq \mathcal {GP}(A)\subseteq \ ^\perp A = \{X\in A\mbox{-}{\rm mod} \ | \ {\rm
Ext}^i_A(X, A) = 0, \ \forall \ i\ge 1\}$. It is clear that
$\mathcal {GP}(A) = A$-mod if and only if $A$ is self-injective. If
$A$ is of finite global dimension then $\mathcal {GP}(A) = \mathcal
{P}(A)$; and if $A$ is {\it a Gorenstein algebra} (i.e., ${\rm
inj.dim}\ _AA< \infty$ and ${\rm inj.dim} \ A_A < \infty$), then \
$\mathcal {GP}(A) = \ ^\perp A$ ([EJ], Corollary 11.5.3). This class
of modules enjoys more stable properties than the usual projective
modules ([AB], where it was called a module of $G$-dimension zero);
it become an important ingredient in the relative homological
algebra ([EJ]) and in the representation theory of algebras (see
e.g. [AR], [B], [GZ], [IKM]); and plays a central role in the Tate
cohomology of algebras (see e.g. [AM] and [Buch]). By [Buch] and
[Hap], the singularity category of Gorenstein algebra $A$ is
triangle equivalent to the stable category of Gorenstein-projective
$A$-modules.

\vskip10pt

On the other hand, the submodule category have been extensively
studied in [RS1] (see also [RS2], [S]). By [KLM] it is also related
to the singularity category (see also [C]). It turns out that the
category of the Gorenstein-projective modules is closly related to
the submodule category, or in general, to the monomorphism category
(see [Z]). The present paper is to explore such an relation in a
more general setting up.

\vskip10pt

Let $\Lambda$ be the path algebra of a finite quiver $Q$ over $A$,
where $A$ is a finite-dimensional algebra over a field $k$. As in
the case of $A=k$, $\Lambda$-modules can be interpreted as
representations of $Q$ over $A$. This interpretation permits us to
introduce the so-called monic representations of $Q$ over $A$. If
$Q= \underset 2 \bullet \longrightarrow \underset 1\bullet$ \ then
they are exactly the objects in the submodule category of $A$ (see
[RS1]); and if $Q= \underset n \bullet \longrightarrow \cdots
\longrightarrow \underset 1\bullet$ \ then the monic representations
of $Q$ over $A$ are exactly the objects in the monomorphism category
of $A$ ([Z]). The main result Theorem \ref{mainthm} of this paper
explicitly determine all the Gorenstein-projective $\m$-modules when
$Q$ is acyclic (i.e., $Q$ has no oriented cycles), via the monic
representations of $Q$ over $A$. We emphasize that here $\m$ is not
necessarily Gorenstein. The proof of Theorem \ref{mainthm} use
induction on $|Q_0|$, and a description of the Gorenstein-projective
modules over the triangular extension of two algebras via bimodules
which is projective in both sides (Theorem \ref{bimodule}). As an
application, we get a characterization of self-injectivity by
claiming that $A$ is self-injective if and only if the
Gorenstein-projective $\m$-modules are exactly the monic
representations of $Q$ over $A$ (Theorem \ref{self}). As another
consequence, if $Q$ has an arrow, then the projective $\m$-modules
coincide with the monic representations of $Q$ over $A$ if and only
$\m$ is hereditary (Theorem \ref{hereditary}).

 \vskip10pt

\section {\bf Monic representations of a quiver over an algebra}

Throughout this section $k$ is a field, $Q$ a finite quiver, and $A$
a finite-dimensional $k$-algebra. We consider the path algebra $AQ$
of $Q$ over $A$, describe its module category, and introduce the
so-called the monic representations of $Q$ over $A$.

\vskip10pt

\subsection{} For the notion of a finite quiver $Q = (Q_0, Q_1, s, e)$ we refer to [ARS] and [R].
We write the conjunction of paths of $Q$ from right to left. Let
$\mathcal P$ be the set of paths of $Q$. Vertex $i$ is a path of
length $0$ and denote it by $e_i$. We define $s(e_i) = i = e(e_i)$.
If $p = \alpha_l\cdots\alpha_1\in\mathcal P$ with $\alpha_i\in Q_1$,
$l\ge 1$, and $e(\alpha_i) = s(\alpha_{i+1})$ for $1\le i\le l-1$,
then we call $l$ the length of $p$ and denote it by $l(p)$, and
define the starting vertex $s(p) = s(\alpha_1)$, and the ending
vertex $e(p) = e(\alpha_l)$. Let $kQ$ be the path algebra of $Q$
over $k$. It is well-known that the category $kQ$-mod of
finite-dimensional $kQ$-modules is equivalent to the category ${\rm
Rep} (Q, k)$ of finite-dimensional representations of $Q$ over $k$.

\vskip10pt

\subsection{} Let $\Lambda = AQ$ be the free left $A$-module with basis $\mathcal P$.
An element of $AQ$ is written as a finite sum
$\sum\limits_{p\in\mathcal P} a_p p$, where $a_p \in A$ and $a_p =
0$ for all but finitely many $p$. Then $\Lambda$ has a $k$-algebra
structure, with multiplication bi-linearly given by $(a_p p)(b_q
q)=(a_p b_q)(pq)$, where $a_pb_q$ is the product in $A$, and $pq$ is
the product in $kQ$. We have isomorphisms  $\Lambda \cong
A\otimes_kkQ \cong kQ\otimes_kA$ of $k$-algebras, and we call
$\Lambda = AQ$ {\it the path algebra of $Q$ over $A$}.  \vskip10pt

For example, if $Q = \underset n \bullet \longrightarrow \cdots
\longrightarrow \underset 1\bullet$ \ then $\Lambda$ is $T_n(A) =
\left(\begin{smallmatrix}
A&A&\cdots&A&A\\
0&A&\cdots&A&A\\
\vdots&\vdots&\vdots&\vdots&\vdots\\
0&0&\cdots&A&A\\ 0&0&\cdots&0&A
\end{smallmatrix}\right)$, the upper triangular matrix algebra of $A$.
In general, if $Q$ is acyclic, then $\Lambda$ is also a kind of
upper triangular matrix algebra over $A$. More precisely, we label
$Q_0$ as $1, \cdots, n$, such that if there is an arrow $\alpha:
j\longrightarrow i$ in $Q_1$ then $j
> i$. Then
$kQ$ is isomorphic to the following matrix algebra over $k$:

$$\left(\begin{smallmatrix}
k&k^{m_{21}}&k^{m_{31}}&\cdots&k^{m_{n1}}\\
0&k&k^{m_{32}}&\cdots&k^{m_{n2}}\\
0&0&k&\cdots &k^{m_{n3}}\\
\vdots&\vdots&\vdots& &\vdots\\
0&0&0&\cdots&k\\
\end{smallmatrix}\right)_{n\times n}, \eqno(2.1)$$
where  $m_{ji}$ is the number of paths from $j$ to $i$, and
$k^{m_{ji}}$ is the direct sum of $m_{ji}$ copies of $k$. It follows
that $\Lambda$ is isomorphic to the following matrix algebra over
$A$:

$$\left(\begin{smallmatrix}
A&A^{m_{21}}&A^{m_{31}}&\cdots&A^{m_{n1}}\\
0&A&A^{m_{32}}&\cdots&A^{m_{n2}}\\
0&0&A&\cdots &A^{m_{n3}}\\
\vdots&\vdots&\vdots& &\vdots\\
0&0&0&\cdots&A\\
\end{smallmatrix}\right)_{n\times n}. \eqno(2.2)$$

\vskip10pt

\subsection{} By definition, {\it a representation $X$ of $Q$ over $A$} is a datum $X
=(X_i, \ X_{\alpha}, \ i\in Q_0, \ \alpha \in Q_1)$, where $X_i$ is
an $A$-module for each  $i\in Q_0$, and $X_{\alpha}: X_{s(\alpha)}
\longrightarrow X_{e(\alpha)}$ is an $A$-map for each $\alpha \in
Q_1$. It is {\it a finite-dimensional representation} if each $X_i$
is finite-dimensional. We will call $X_i$ {\it the $i$-th branch} of
$X$. A morphism $f$ from representation $X$ to representation $Y$ is
a datum $(f_i, \ i\in Q_0)$, where $f_i: X_i \longrightarrow Y_i$ is
an $A$-map for each $i\in Q_0$, such that for each arrow $\alpha:
j\longrightarrow i$ the following diagram

\[\xymatrix {X_j\ar[r]^{f_j}\ar[d]^{X_\alpha} & Y_j\ar[d]^{Y_\alpha} \\
X_i\ar[r]^{f_i} & Y_i}\eqno(2.3)\] commutes. If $p =
\alpha_l\cdots\alpha_1\in\mathcal P$ with $\alpha_i\in Q_1$,  $l\ge
1$, and $e(\alpha_i) = s(\alpha_{i+1})$ for $1\le i\le l-1$, then we
put $X_p$ to be the $A$-map $X_{\alpha_l}\cdots X_{\alpha_1}$.
Denote by ${\rm Rep}(Q, A)$ the category of finite-dimensional
representations of $Q$ over $A$. A morphism $f = (f_i, \ i\in Q_0)$
in ${\rm Rep}(Q, A)$ is a monomorphism (resp., an epimorphism, an
isomorphism) if and only if for each $i\in Q_0$, $f_i$ is injective
(resp., surjective, an isomorphism).

\vskip10pt

\begin{lem} \label{rep} \ Let $\Lambda$ be the path algebra of
$Q$ over $A$. Then we have an equivalence $\Lambda\mbox{-}{\rm
mod}\cong {\rm Rep}(Q, A)$ of categories, here $\Lambda\mbox{-}{\rm
mod}$ is the category of finite-dimensional $\m$-modules.\end{lem}

We omit the details of the proof of Lemma \ref{rep}, which is
similar to the case of $A=k$ (see Theorem 1.5 of [ARS], p.57; or
[R], p.44). In the following we will identify a $\Lambda$-module
with a representation of $Q$ over $A$, which is always assumed to be
finite-dimensional. Under this identification, a $\m$-module $X$ is
a representation $(X_i, \ X_\alpha, \ i\in Q_0, \ \alpha\in Q_1)$ of
$Q$ over $A$, where $X_i = (1e_i)X$,  \ $1$ is the identity of $A$,
and the $A$-action on $X_i$ is given by $a (1e_i)x = (ae_i)x =
(1e_i)(ae_i)x, \ \forall \ x\in X, \ \forall \ a\in A$; and
$X_\alpha: X_{s(\alpha)}\longrightarrow X_{e(\alpha)}$ is the
$A$-map given by the left action by $1\alpha\in \m$.  On the other
hand, a representation $(X_i, \ X_\alpha, \ i\in Q_0, \ \alpha\in
Q_1)$ of $Q$ over $A$ is a $\m$-module $X = \bigoplus\limits_{i\in
Q_0}X_i,$ with the $\m$-action on $X$ given by
$$(ap)(x_i) = \begin{cases} 0, & \mbox{if} \ s(p)\ne i;
\\ ax_i, & \mbox{if} \ s(p) = i, \ l(p) = 0; \\ aX_p(x_i)\in X_{e(p)}, & \mbox{if} \ s(p) =
i, \ l(p) \ge 1,
\end{cases}  \ \ \ \ \forall \ a\in A, \ \forall \ p\in\mathcal P, \ \forall \
x_i\in X_i.$$ An indecomposable projective $\m$-module is of the
form $L\otimes_k P(i)$, where $P(i)$ is the indecomposable
$kQ$-module at vertex $i$, and $L$ is an indecomposable projective
$A$-module. In particular, each branch of a projective $\m$-module
is a projective $A$-module.

\vskip10pt

Let $f: X\longrightarrow Y$ be a morphism in ${\rm Rep}(Q, A)$. Then
${\rm Ker} f$ and $\cok f$ can be explicitly written out. For
example, $\cok f = (\cok f_i, \ \widetilde{Y_{\alpha}}, \ i\in Q_0,
\ \alpha\in Q_1)$, where for each arrow $\alpha: j \longrightarrow
i$, $\widetilde{Y_{\alpha}}: \cok f_j \longrightarrow \cok f_i$ is
the $A$-map induced by $Y_\alpha$ (see $(2.3)$). A sequence of
morphisms $0\longrightarrow X\stackrel {f}\longrightarrow Y\stackrel
{g}\longrightarrow Z \longrightarrow 0$ in ${\rm Rep}(Q, A)$ is
exact if and only if $0\longrightarrow X_i\stackrel
{f_i}\longrightarrow Y_i\stackrel {g_i}\longrightarrow Z_i
\longrightarrow 0$ is exact in $A$-mod for each $i\in Q_0$.

\vskip10pt

In the following, if $Q_0$ is labeled as $1, \cdots, n$, then we
also write a representation $X$ of $Q$ over $A$ as
$\left(\begin{smallmatrix}
X_1\\
\vdots\\
X_n
\end{smallmatrix}\right)_{(X_\alpha, \ \alpha\in Q_1)}$.

\vskip10pt

\subsection{} The following is a central notion in this paper.

\vskip10pt

\begin{defn} \label{maindef} \ A representation $X = (X_i,
X_{\alpha}, i\in Q_0, \alpha\in Q_1)$ of $Q$ over $A$ is {\it a
monic representation}, or a monic $\Lambda$-module, if for each
$i\in Q_0$ the following $A$-map
$$(X_{\alpha})_{\alpha\in Q_1, \ e(\alpha) = i}: \ \bigoplus\limits_{\begin
{smallmatrix} \alpha\in Q_1\\ e(\alpha) = i \end{smallmatrix}}
X_{s(\alpha)} \longrightarrow X_i$$ is injective,  or equivalently,
the following two conditions are satisfied

\vskip5pt

$(m1)$  \ For each $\alpha\in Q_1$, $X_{\alpha}: X_{s(\alpha)}
\longrightarrow X_{e(\alpha)}$ is an injective map; and

\vskip5pt

$(m2)$  \ For each $i\in Q_0$,  there holds  $\sum\limits_{\begin
{smallmatrix} \alpha\in Q_1\\ e(\alpha) = i \end{smallmatrix}}\Ima
X_\alpha = \bigoplus\limits_{\begin {smallmatrix} \alpha\in Q_1\\
e(\alpha) = i \end{smallmatrix}}\Ima X_\alpha.$
\end{defn}

\vskip10pt

Denote by ${\rm Mon}(Q, A)$ the full subcategory of ${\rm Rep}(Q,
A)$ consisting of the monic representations of $Q$ over $A$. In
particular, if $A=k$, then we have ${\rm Mon}(Q, k)\subseteq {\rm
Rep}(Q, k)$.

\vskip10pt

For example, if $Q = \underset n \bullet \longrightarrow \cdots
\longrightarrow \underset 1\bullet$ \ \ then a representation $X$ of
$Q$ over $A$ is simply written as $X = (X_i, \ \phi_i)$, where each
$\phi_i: X_{i+1} \longrightarrow X_i$ is an $A$-map, $1\le i\le
n-1$. In this case $X$ is a monic representation, or a monic
$T_n(A)$-module, exactly means that each $\phi_i $ is injective,
$1\le i\le n-1$. This kind of monic $T_n(A)$-modules have arisen
from different questions and in different terminologies, for
examples in [RS1], [RS2], [S], [LZ], [KLM], [C], [Z], [IKM].

\vskip10pt

\subsection{}

There is another similar but different notion. Let $A=kQ/I$ be a
finite-dimensional $k$-algebra, where $I$ is an admissible ideal of
$kQ$. An $I$-bounded representations of $Q$ over $k$ is a datum $X
=(X_i, \ X_{\alpha}, \ i\in Q_0, \ \alpha \in Q_1)$, where $X_i$ is
a $k$-space for each $i\in Q_0$, and $X_{\alpha}: X_{s(\alpha)}
\longrightarrow X_{e(\alpha)}$ is a $k$-linear map for each $\alpha
\in Q_1$, such that $\sum\limits_{p\in\mathcal P}c_pX_p = 0$ for
each element $\sum\limits_{p\in\mathcal P}c_pp\in I$, where $l(p)\ge
2$ and $c_p\in k$. An $I$-bounded representation $X = (X_i,
X_{\alpha}, i\in Q_0, \alpha\in Q_1)$ of $Q$ over $k$ is {\it a
monic representation}, if for each $i\in Q_0$ the following
$k$-linear map
$$(X_{\alpha})_{\alpha\in Q_1, \ e(\alpha) = i}: \ \bigoplus\limits_{\begin
{smallmatrix} \alpha\in Q_1\\ e(\alpha) = i \end{smallmatrix}}
X_{s(\alpha)} \longrightarrow X_i$$ is injective. Let ${\rm Rep}(Q,
I, k)$ be the category of finite-dimensional $I$-bounded
representations of $Q$ over $k$. It is well-known that there is an
equivalence $A$-mod$\cong {\rm Rep}(Q, I, k)$ of categories (see
Proposition 1.7 in [ARS], p.60; or [R], p.45). Let ${\rm Mon}(Q, I,
k)$ denote the full subcategory of ${\rm Rep}(Q, I, k)$ of
$I$-bounded monic representations $Q$ over $k$. Thus ${\rm Mon}(Q,
0, k) = {\rm Mon}(Q, k).$

\vskip10pt

\begin{prop} \label{Ahereditary} \  Let $A=kQ/I$ be a
finite-dimensional $k$-algebra, where $I$ is an admissible ideal of
$kQ$. Then $\mathcal{P}(A) \subseteq {\rm Mon}(Q, I, k)$ if and only
if $A$ is hereditary.
\end{prop}

\noindent{\bf Proof.} \  If $A$ is hereditary then $I = 0$. It is
clear $\mathcal{P}(kQ) \subseteq {\rm Mon}(Q, 0, k)$.

Conversely, if $I\ne 0$, then take an element $\sum\limits_{p\in
\mathcal P}c_pp\in I$ with $l(p) \ge 2$ and $c_p\in k$. Assume that
all the paths $p$ with $c_p\ne 0$ have the same starting vertex $j$
and the same ending vertex $i$. Consider the projective $A$-module
$P(j) = Ae_j$. As an $I$-bounded representation of $Q$ over $k$ we
have $P(j) = (e_tkQe_j,
 t\in Q_0, \ f_\alpha, \ \alpha\in Q_1)$. Let $\alpha_1, \cdots,
\alpha_m$ be all the arrows of $Q$ ending at $i$. We claim that
$$(f_{\alpha_v})_{1\le v\le m}: \ \bigoplus\limits_{1\le v\le m}
e_{s(\alpha_v)}kQe_j \longrightarrow e_ikQe_j$$ is not injective,
where $f_{\alpha_v}$ is the $k$-linear map given by the left
multiplication by $\alpha_v$. Since each path from $j$ to $i$ must
go through some $\alpha_v$, and $\sum\limits_{p\in \mathcal P}c_pf_p
= 0$, it follows that $\sum\limits_{1\le v\le m} {\rm
dim}_k(e_{s(\alpha_v)}kQe_j) > {\rm dim}_k (e_ikQe_j).$ This
justifies the claim, i.e., $P(j)\notin {\rm Mon}(Q, I, k)$. $\s$

\vskip10pt

Now, let $\m = A\otimes_kkQ$ be the path algebra of $Q$ over $A$.
Assume that $\m$ is of the form $\m = kQ'/I'$, where $Q'$ is a
finite quiver and $I'$ is an admissible ideal of $kQ'$. We emphasize
that in general ${\rm Mon}(Q, A) \ne {\rm Mon}(Q', I', k)$. In fact,
we will see in Theorem \ref{mainthm} that $\mathcal{P}(\m) \subseteq
{\rm Mon}(Q, A)$ is always true; but in general $\mathcal{P}(\m)
\subseteq {\rm Mon}(Q', I', k)$ is not true,  as Proposition
\ref{Ahereditary} shows. This is the reason why we do not use the
notation ${\rm Mon}(\m)$.

\vskip10pt

\section {\bf Algebras given by bimodules}

\vskip10pt

\subsection{} \ Let $A$ and $B$ be rings, and $M$ an $A$-$B$-bimodule. Consider the
upper triangular matrix ring $\m=\left(\begin{smallmatrix}
A&M\\
0&B
\end{smallmatrix}\right)$,
where the addition and the multiplication are given by the ones of
matrices. We assume that $\m$ is an Artin algebra ([ARS], p.72), and
only consider finitely generated $\m$-modules. A $\m$-module can be
identified with a tripe $\left(\begin{smallmatrix}
X\\
Y
\end{smallmatrix}\right)_\phi$, or simply  $\left(\begin{smallmatrix}
X\\
Y
\end{smallmatrix}\right)$
if $\phi$ is clear, where $X\in A$-mod, $Y\in B$-mod, and $\phi:
M\otimes_{B} Y\longrightarrow X$ is an $A$-map. A $\m$-map
$\left(\begin{smallmatrix}
X \\
   Y
 \end{smallmatrix}\right)_{\phi}\longrightarrow \left(\begin{smallmatrix}
   X'  \\
   Y'
 \end{smallmatrix}\right)_{\phi'}$
 can be identified with a pair
 $\left(\begin{smallmatrix}
   f  \\
   g
 \end{smallmatrix}\right)$,
where $f\in {\Hom}_A (X, X'),$ \ $g\in \Hom_B (Y, Y')$, such that
the diagram \ $$\xymatrix{
M\otimes_B Y\ar[d]_{{\rm id}\otimes g} \ar[r]^-\phi & X \ar[d]_f \\
M\otimes_B Y' \ar[r]^-{\phi'} & X'}$$ commutes. A sequence of
$\Lambda$-maps $0\longrightarrow \left(\begin{smallmatrix}
X_1 \\
   Y_1
 \end{smallmatrix}\right)_{\phi_1}\stackrel {\left(\begin{smallmatrix}
   f_1  \\
   g_1
 \end{smallmatrix}\right)}\longrightarrow \left(\begin{smallmatrix}
   X_2  \\
   Y_2
 \end{smallmatrix}\right)_{\phi_2} \stackrel {\left(\begin{smallmatrix}
   f_2  \\
   g_2
 \end{smallmatrix}\right)}\longrightarrow \left(\begin{smallmatrix}
   X_3  \\
   Y_3
 \end{smallmatrix}\right)_{\phi_3} \longrightarrow 0$ is exact if and only if
$0\longrightarrow X_1\stackrel {f_1}\longrightarrow X_2\stackrel
{f_2}\longrightarrow X_3 \longrightarrow 0$ is an exact sequence of
$A$-maps, and $0\longrightarrow Y_1\stackrel {g_1}\longrightarrow
Y_2\stackrel {g_2}\longrightarrow Y_3 \longrightarrow 0$ is an exact
sequence of $B$-maps. Indecomposable projective $\m$-modules are
exactly $\left(\begin{smallmatrix}
   P \\
   0
\end{smallmatrix}\right)
$ and $\left( \begin{smallmatrix}
   M\otimes_B Q  \\
   Q
\end{smallmatrix}\right)_{\rm id}$,
where $P$ runs over indecomposable projective $A$-modules, and $Q$
runs over indecomposable projective $B$-modules.

\vskip10pt

\subsection{} The following result describes the Gorenstein-projective
$\Lambda$-modules, if $_AM$ and $M_B$ are projective modules. We
emphasize that here $\m$ is not assumed to be Gorenstein (see
Corollary 3.3 of [XZ] for the similar result under the assumption
that $\m$ is Gorenstein; and the proof there in [XZ] can not be
generalized to the non-Gorenstein case).

\vskip10pt

\begin{thm} \label{bimodule} \ Let $\m = \left(\begin{smallmatrix}
A&M\\
0&B
\end{smallmatrix}\right)$ be an Artin algebra, $M$
an $A$-$B$-bimodule such that $_AM$ and $M_B$ are projective
modules.  Then
$\left(\begin{smallmatrix}X\\Y\end{smallmatrix}\right)_\phi\in
\mathcal {GP}(\m)$ if and only if $\phi: M\otimes_B Y\longrightarrow
X$ is injective, $\cok\phi\in \mathcal {GP}(A)$, and $Y\in \mathcal
{GP}(B)$. In this case, $X\in \mathcal {GP}(A)$ if and only if
$M\otimes_B Y\in \mathcal {GP}(A)$.
\end{thm}

\noindent{\bf Proof.} \ The last assertion is easy, since in this
case $0 \longrightarrow M\otimes_B Y \stackrel \phi \longrightarrow
X \longrightarrow \cok\phi \longrightarrow 0$ is exact, and
$\mathcal {GP}(A)$ is  closed under extensions and the kernels of
epimorphisms (see e.g. [Hol]).

\vskip5pt

We first prove the ``if" part. Assume that $\phi: M\otimes_B
Y\longrightarrow X$ is injective, $\cok\phi\in \mathcal {GP}(A)$,
and $Y\in \mathcal {GP}(B)$. Then we have a complete $B$-projective
resolution
$$Q^\bullet = \cdots\To Q^{-1}\To Q^0\stackrel{d'^0}\To Q^1\To \cdots \eqno(3.1)$$
with $Y = \Ker d'^0$, and a complete $A$-projective resolution
$$P^\bullet = \cdots\To P^{-1}\To P^0\stackrel{d^0}\To P^1\To \cdots \eqno(3.2)$$
with $\cok \phi = \Ker d^0$. Since $M_B$ is projective, we get the
following exact sequences of $A$-modules

\begin{align*} & 0 \longrightarrow M\otimes_B Y\longrightarrow
M\otimes_B  Q^0\longrightarrow M\otimes_B Q^1 \longrightarrow \cdots
\\ & 0 \longrightarrow \cok\phi\longrightarrow \ \ \ P^0 \ \ \ \
\longrightarrow \ \ \   P^1 \ \ \ \ \longrightarrow  \cdots.
\end{align*}
Since $_AM$ is projective, $M\otimes_B Q^i$ is a projective
$A$-module for each $i\ge 0$. Note that projective $A$-modules are
injective objects in $\mathcal {GP}(A)$, it follows from the exact
sequence $0\longrightarrow M\otimes_B Y\longrightarrow X
\longrightarrow \cok\phi\longrightarrow 0$ and a version of
Horseshoe Lemma that there is an exact sequence of $A$-modules
$$0\longrightarrow X \longrightarrow P^0\oplus (M\otimes_B
Q^0) \stackrel {\partial^0} \longrightarrow P^1\oplus (M\otimes_B
Q^1)\longrightarrow \cdots \eqno(3.3)$$ with $\partial^i =
\left(\begin{smallmatrix}
d^i & 0\\
\sigma^i & {\rm id}\otimes_B d'^i
\end{smallmatrix}\right),$ \ $\sigma^i: P^i \longrightarrow M\otimes_BQ^i, \ \forall \ i\in \Bbb Z,$
such that the following diagram

\[\xymatrix {\ \ \ \ 0 \ar[r] & M\otimes_BY\ar[d]^-{\phi}\ar[r]& M\otimes_B
Q^0\ar[d]_{\binom{0}{\rm id}}\ar[r]^{{\rm id}\otimes_B d'^0}&
M\otimes_B Q^1 \ar[d]^-{\binom{0}{\rm id}}\ar[r] & \cdots \ \ \ \ \
\ \ \ \ \ \ \ \ \ \ \ \ \ \ \ \ \ \ \    \\ \ \ \ \ \ \ \ \ 0 \ar[r]
& X\ar[r] & P^0\oplus (M\otimes_B Q^0)\ar[r]^{\partial^0}& P^1\oplus
(M\otimes_B Q^1)\ar[r] & \cdots \ \ \ \ \ \ \ \ \ \ \ (3.4)}\]
commutes. By the same argument we get the following commutative
diagram with exact rows

\[\xymatrix {\ \ \ \cdots \ar[r] & M\otimes_BQ^{-2} \ar[d]^-{\binom{0}{\rm id}} \ar[r]^{{\rm id}\otimes_Bd'^{-2}}
 &  M\otimes_B
Q^{-1} \ar[d]_{\binom{0}{\rm id}} \ar[r] & M\otimes_B Y \ar[d]^-
\phi\ar[r] & 0 \ \ \ \ \ \ \ \ \   \\ \ \ \ \cdots  \ar[r] &
P^{-2}\oplus (M\otimes_B Q^{-2}) \ar[r]^{\partial^{-2}} &
P^{-1}\oplus (M\otimes_B Q^{-1}) \ar[r] & X \ar[r] & 0. \ \ \ \ \ \
\ \ \ (3.5)}
\]
Putting $(3.4)$ and $(3.5)$ together we then get the following exact
sequence of projective $\m$-modules
$$L^\bullet = \cdots \longrightarrow
\left(\begin{smallmatrix}
P^{-1}\oplus (M\otimes_{B}Q^{-1})\\
Q^{-1}
\end{smallmatrix}\right){\longrightarrow}
\left(\begin{smallmatrix}
P^0\oplus (M\otimes_{B}Q^{0})\\
Q^{0}
\end{smallmatrix}\right)_{\binom{0}{\rm id}}\stackrel {\left(\begin{smallmatrix}
\partial^0\\
d'^0
\end{smallmatrix}\right)}{\longrightarrow} \left(\begin{smallmatrix}
P^1\oplus (M\otimes_{B}Q^{1})\\
Q^1
\end{smallmatrix}\right)
{\longrightarrow} \cdots \eqno (3.6)$$ with $\Ker
{\left(\begin{smallmatrix}
\partial^0\\
d'^0
\end{smallmatrix}\right)}
=\left(\begin{smallmatrix}X\\Y\end{smallmatrix}\right)_\phi$.

For each projective $A$-module $P$, $\Hom_{\m}(L^\bullet,
\left(\begin{smallmatrix}
P\\
0
\end{smallmatrix}\right)) \cong \Hom_{A}(P^\bullet, P)$ is exact, since $P^\bullet$ is a complete projective resolution.
For each projective $B$-module $Q$, since $Q^\bullet$ is a complete
projective resolution, $\Hom_{B}(Q^\bullet, Q)$ is exact. Since
$M\otimes_B Q$ is projective, $\Hom_{A}(P^\bullet, M\otimes_BQ)$ is
exact. Note that $$\Hom_{\m}(L^\bullet, \left(\begin{smallmatrix}
M\otimes_B Q\\
Q
\end{smallmatrix}\right)) \cong \Hom_{A}(P^\bullet,
M\otimes_BQ)\oplus \Hom_{B}(Q^\bullet, Q)$$ (here the direct sum
only means that each term of the complex at the left hand side is a
direct sum of terms of complexes at the right hand side, i.e., it
does not mean a direct sum of complexes. In fact, the complex at the
right hand side has differentials $\left(\begin{smallmatrix}
\Hom_A(d^i, M\otimes_BQ) & \Hom_A(\sigma^i, M\otimes_BQ)\\ 0 &
\Hom_B(d'^i, Q)
\end{smallmatrix}\right)$). By the canonical exact sequence of complexes
$$0\longrightarrow \Hom_{A}(P^\bullet,
M\otimes_BQ) \stackrel{\left(\begin{smallmatrix}
{\rm id}\\
0
\end{smallmatrix}\right)} \longrightarrow \Hom_{\m}(L^\bullet, \left(\begin{smallmatrix}
M\otimes_B Q\\
Q
\end{smallmatrix}\right))
\stackrel{\left(\begin{smallmatrix} 0 & {\rm id}
\end{smallmatrix}\right)}\longrightarrow \Hom_{B}(Q^\bullet, Q) \longrightarrow 0$$ and the fundamental
theorem of homological algebra we see that $\Hom_{\m}(L^\bullet,
\left(\begin{smallmatrix}
M\otimes_B Q\\
Q
\end{smallmatrix}\right))$ is also
exact. Therefore we conclude that $L^\bullet$ is a complete
$\m$-projective resolution, and hence
$\left(\begin{smallmatrix}X\\Y\end{smallmatrix}\right)_\phi$ is a
Gorenstein-projective $\Lambda$-module.

\vskip10pt

Conversely, assume that
$\left(\begin{smallmatrix}X\\Y\end{smallmatrix}\right)_\phi\in
\mathcal {GP}(\m)$. Then there is a complete $\m$-projective
resolution $(3.6)$ with $\Ker {\left(\begin{smallmatrix}
\partial^0\\
d'^0
\end{smallmatrix}\right)}
=\left(\begin{smallmatrix}X\\Y\end{smallmatrix}\right)_\phi$. Then
we get an exact sequence $(3.1)$ of projective $B$-modules with
$\Ker d'^0 = Y$,  and the following exact sequence

$$V^\bullet = \cdots \longrightarrow P^{-1}\oplus (M\otimes_B Q^{-1})
\longrightarrow P^0\oplus (M\otimes_B Q^0) \stackrel {\partial^0}
\longrightarrow P^1\oplus (M\otimes_B Q^1)\longrightarrow \cdots
\eqno(3.7)$$ of projective $A$-modules with $\Ker \partial^0 = X$.
Since $M_B$ is projective, it follows that $M\otimes_BQ^\bullet$ is
exact. Since ${\left(\begin{smallmatrix}
\partial^i\\
d'^i
\end{smallmatrix}\right)}$ is a $\m$-map, by $(3.6)$ we know that $\partial^i$ is of
the form $\partial^i = \left(\begin{smallmatrix}
d^i & 0\\
\sigma^i & {\rm id}\otimes_B d'^i
\end{smallmatrix}\right),$ where $\sigma^i: P^i \longrightarrow M\otimes_BQ^i, \ \forall \ i\in \Bbb Z,$ and
$$P^\bullet = \cdots\To P^{-1}\To P^0\stackrel{d^0}\To P^1\To \cdots$$
is a complex. By the canonical exact sequence of complexes
$$0\longrightarrow M\otimes_BQ^\bullet \stackrel{\left(\begin{smallmatrix}
0\\
{\rm id}
\end{smallmatrix}\right)} \longrightarrow V^\bullet
\stackrel{\left(\begin{smallmatrix} {\rm id}, & 0
\end{smallmatrix}\right)}\longrightarrow P^\bullet \longrightarrow 0$$ and the fundamental
theorem of homological algebra we see that $P^\bullet$ is also
exact.

From $(3.6)$ we have the following commutative diagram with exact
rows and columns:

$$
\CD && && 0 && 0\cr &&&& @VVV @VVV \cr 0@>>> M\otimes_{B}Y @>>>
M\otimes_{B}Q^{0} @>>> M\otimes_{B}Q^{1}@>>> \cdots \cr & & @VV\phi
V @VV\left(\begin{smallmatrix}
0\\
{\rm id}
\end{smallmatrix}\right) V   @VV\left(\begin{smallmatrix}
0\\
{\rm id}
\end{smallmatrix}\right)V   \cr
  0 @>>> X  @>>> P^{0}\oplus (M\otimes_{B}Q^{0})@>>> P^{1}\oplus (M\otimes_{B}Q_{1})@>>> \cdots \cr
 & &
@VV V     @VV({\rm id},0) V   @VV({\rm id},0)V   \cr
  0 @>>> \cok \phi @>>> P^{0}@>d^0>> P^{1}@>>> \cdots
\cr& & @VVV     @VV V   @VVV   \cr
  && 0 & & 0 && 0&
\endCD
$$

\vskip5pt

\noindent Thus $\phi: M\otimes_BY\To X$ is injective and $\Ker
{d}^0\cong \cok \phi$. For each projective $A$-module $P$, since
$\Hom_{\m}(L^\bullet, \left(\begin{smallmatrix}
P\\
0
\end{smallmatrix}\right)) \cong \Hom_{A}(P^\bullet, P)$ and
$L^\bullet$ is a complete projective resolution, it follows that
$P^\bullet$ is a complete projective resolution, and hence
$\cok\phi$ is a Gorenstein-projective $A$-module.

For each projective $B$-module $Q$, since $P^\bullet$ is a complete
projective resolution, it follows that $\Hom_{A}(P^\bullet,
M\otimes_BQ)$ is exact. Since $L^\bullet$ is a complete projective
resolution, it follows that
$$\Hom_{\m}(L^\bullet, \left(\begin{smallmatrix}
M\otimes_B Q\\
Q
\end{smallmatrix}\right)) \cong \Hom_{A}(P^\bullet,
M\otimes_BQ)\oplus \Hom_{B}(Q^\bullet, Q)$$ is exact (again, here
the direct sum does not mean a direct sum of complexes). By the same
argument we know that $\Hom_{B}(Q^\bullet, Q)$ is exact. It follows
that $Y$ is a Gorenstein-projective $B$-module. This completes the
proof. $\s$

\vskip5pt

We remark that if $\m$ is Gorenstein, then in Theorem \ref{bimodule}
$\left(\begin{smallmatrix}X\\Y\end{smallmatrix}\right)_\phi\in
\mathcal {GP}(\m)$ implies $X\in \mathcal {GP}(A)$ (see [XZ],
Corollary 3.3).

\section {\bf Main result}

\vskip10pt

\subsection{} \ The aim of this section is to prove the following
characterization of Gorestein-projective $\m$-modules, where $\m$ is
the path algebra of a finite acyclic quiver over a
finite-dimensional algebra. We emphasize that here $\m$ is not
assumed to be Gorenstein.

\vskip5pt

\begin{thm} \label{mainthm} \ Let $Q$ be a finite acyclic quiver, and $A$ a
finite-dimensional algebra over a field $k$. Let $\m =
A\otimes_kkQ$,  and $ X = (X_i, \ X_\alpha, \ i\in Q_0, \ \alpha\in
Q_1)$ be a $\m$-module. Then $X\in\mathcal {GP}(\m)$ if and only if
$X\in {\rm Mon}(Q, A)$ and $X$ satisfies the following condition
$(G)$, where

\vskip5pt

$(G)$ \ For each $i \in Q_0$,  $X_i\in\mathcal {GP}(A)$, and the quotient ${X_i}/ ({\bigoplus\limits_{\begin {smallmatrix} \alpha\in Q_1\\
e(\alpha) = i \end{smallmatrix}}\Ima X_\alpha})\in \mathcal
{GP}(A)$.
\end{thm}

\vskip10pt

\begin{exm}  $(i)$ \  Taking $Q = \underset n \bullet \longrightarrow \cdots \longrightarrow
 \underset 1\bullet$ in Theorem
\ref{mainthm} we get: a $T_n(A)$-module $X =  (X_i, \phi_i)$ is
Gorenstein-projective if and only if each $\phi_i$ is injective,
each $X_i$ is a Gorenstein-projective $A$-module, and each
$\cok\phi_i$ is a Gorenstein-projective $A$-module. Under the
assumption that $A$ is Gorenstein, this result was obtained in
Corollary 4.1 of [Z].

\vskip10pt

$(ii)$ \ Let $\Lambda$ be the $k$-algebra given by quiver
$\xymatrix{\underset 3 {\bullet}
\ar@(ul,ur)[]^{\lambda_3}\ar[r]^\beta &\underset 1 \bullet
\ar@(ul,ur)[]^{\lambda_1} & \underset
2\bullet\ar@(ul,ur)^{\lambda_2}\ar[l]_\alpha}$ with relations
$\lambda_1^2, \ \lambda_2^2, \ \lambda_3^2, \ \alpha\lambda_2
-\lambda_1\alpha, \ \beta\lambda_3 -\lambda_1\beta$. Then $\Lambda =
A\otimes_k kQ = \left (\begin {smallmatrix} A & A& A\\ 0 & A & 0 \\
0& 0& A
\end{smallmatrix}\right)$, where $Q$ is the quiver $\underset 3
\bullet \longrightarrow \underset 1 \bullet \longleftarrow \underset
2 \bullet$, \  and $A = k[x]/\langle x^2\rangle$. Let $k$ be the
simple $A$-module, and $\sigma: k\hookrightarrow A$ the inclusion.
Then by Theorem \ref{mainthm}
$$X = (X_1
= A\oplus k, \ X_2 = k, \ X_3 = k, \ X_\alpha = \binom {0}{\rm id},
\ X_\beta = \binom{\sigma}{{\rm id}})\in\mathcal {GP}(\m);$$ while
$$Y = (Y_1 = A, \ Y_2 = k, \ Y_3 = k, \ Y_\alpha =
\sigma = Y_\beta)\notin\mathcal {GP}(\m).$$
\end{exm}

\vskip10pt

\subsection{} \ Theorem \ref{mainthm} will be proved by
using Theorem \ref{bimodule} and induction on $|Q_0|$, the number of
vertices of $Q$.

\vskip10pt

We label $Q_0$ as $1, \cdots, n$, such that if there is an arrow
$\alpha: j\longrightarrow i$ in $Q_1$, then $j
> i$. Thus $n$ is a source of $Q$. Denote  by $Q'$ the
quiver obtained from $Q$ by deleting vertex $n$, and by $\Lambda' =
A\otimes_kkQ'$ the path algebra of $Q'$ over $A$. Let $P(n)$ be the
indecomposable projective (left) $kQ$-module at vertex $n$. Put $P =
A\otimes_k{\rm rad} P(n)$. Clearly $P$ is a $\m'$-$A$-bimodule and
$\Lambda =\left(\begin{smallmatrix}
\Lambda' &P\\
0&A
\end{smallmatrix}\right)$. See $(2.2)$.

\vskip10pt

Since $kQ$ is hereditary, ${\rm rad} P(n)$ is a  projective
$kQ'$-module, and hence $P = A\otimes _k{\rm rad} P(n)$ is a (left)
projective $\Lambda'$-module, and a (right) projective $A$-module
(since as a right $A$-module, $P$ is a direct sum of copies of
$A_A$). This allows us to apply Theorem \ref{bimodule}. For this, we
write a $\Lambda$-module $X= (X_i, \ X_\alpha, \ i\in Q_0, \
\alpha\in Q_1)$ as $X = \left(\begin{smallmatrix}
X'\\
X_n
\end{smallmatrix}\right)_{\phi}$,
where $X'= (X_i, \ X_\alpha, \ i\in Q'_0, \ \alpha\in Q'_1)$ is a
$\Lambda'$-module, and  $\phi: P\otimes_A X_n\longrightarrow X'$ is
a $\Lambda'$-map, whose explicit expression will be given in the
proof of Lemma \ref{mono} below.

\vskip10pt

We will keep all these notations of $Q'$, $\Lambda'$, $P(n), \ P$,
$X'$ and $\phi$, throughout this section.

\subsection{}
\ By a direct translation from Theorem \ref{bimodule} in this
special case, we have

\vskip10pt

\begin{lem} \label{mainlem} \ Let $X= \left(\begin{smallmatrix}
X'\\
X_n
\end{smallmatrix}\right)_{\phi}$ be a $\Lambda$-module.
Then $X\in\mathcal {GP}(\m)$  if and only if $X$ satisfies the
following conditions:

\vskip5pt

$(i)$ \ $X_n\in\mathcal {GP}(A)$;

\vskip5pt

$(ii)$ \ $\phi: P\otimes_A X_n \longrightarrow X'$ is injective;

\vskip5pt

$(iii)$ \ $\cok\phi\in\mathcal {GP}(\m')$.

\vskip5pt

\end{lem}

\vskip10pt

For each $i \in Q'_0$, denote by $\ha(n\to i)$ the set of the arrows
from $n$ to $i$; and by $\hp(n\to i)$ the set of paths from $n$ to
$i$. For an integer $m\ge 0$ and a module $M$, let $M^m$ denote the
direct sum of $m$ copies of $M$.

\vskip5pt

\begin{lem}\label{mono} \ Let $ X = (X_i, \ X_\alpha, \ i\in Q_0, \
\alpha\in Q_1)$ be a $\m$-module. If $\ X_\beta$ is injective for
each $\beta \in Q'_1$, then \ $\phi: P\otimes_A X_n \longrightarrow
X'$ is injective if and only if $\ X_\alpha$ is injective, $\forall
\ \alpha \in Q_1$, and $\sum\limits_{p \in \hp(n \to i)}\Ima X_p
=\bigoplus\limits_{p \in \hp(n \to i)}\Ima X_p,  \ \forall \ i\in
Q_0'.$
\end{lem}

\noindent {\bf Proof.} \ For $i \in Q'_0$, set $m_i= |\mathcal{P}(n
\to i)|$. As a $kQ'$-module, ${\rm rad} P(n)$  can be written as
$\left(\begin{smallmatrix}
k^{m_1}\\
\vdots\\
k^{m_{n-1}}\\
\end{smallmatrix}\right)$ (please see (2.1) and 4.2), hence we have isomorphisms of $\m'$-modules
$$P\otimes_A X_n \cong ({\rm rad}
P(n)\otimes_kA)\otimes_A X_n \cong {\rm rad} P(n)\otimes_k X_n \cong \left(\begin{smallmatrix} X_n^{m_1}\\
\vdots \\
X_n^{m_{n-1}}\\
\end{smallmatrix}\right).$$ Let $\hp(n\to i) = \{p_1, \cdots,p_{m_i}\}$.
Then $\phi$ is of the form
$$\left(\begin{smallmatrix}
\phi_1\\
\vdots\\
\phi_{n-1}\\
\end{smallmatrix}\right): P\otimes_A X_n \cong \left(\begin{smallmatrix}
X_n^{m_1}\\
\vdots\\
X_n^{m_{n-1}}\\
\end{smallmatrix}\right)
\longrightarrow  \left(\begin{smallmatrix}
X_1\\
\vdots\\
X_{n-1}\\
\end{smallmatrix}\right)
,$$ where $\phi_i=(X_{p_1}, \ \cdots,  \ X_{p_{m_i}}): X_n^{m_i}
\longrightarrow X_i$ (for the meaning of $X_{p_i}$ please see 2.3).
So $\phi$ is injective if and only if $\phi_i$ is injective for each
$i\in Q'_0$; if and only if $\sum\limits_{p \in \mathcal{P}(n \to
i)}\Ima X_p =\bigoplus\limits_{p \in \hp(n \to i)}\Ima X_p $, and $\
X_p$ is injective, $\forall \ p \in \mathcal P(n \to i)$,  from
which and the assumption the assertion follows. $\s$

\vskip10pt

\begin{lem} \label{mono2} Let $X=\left(\begin{smallmatrix}
X'\\
X_n
\end{smallmatrix}\right)_{\phi}$ be a monic $\Lambda$-module. Then

\vskip5pt

$(1)$ \ For each $i\in Q_0'$ there holds $\sum\limits_{p\in
\hp(n\rightarrow i)}\Ima X_p = \bigoplus\limits_{p\in
\hp(n\rightarrow i)}\Ima X_p$;

\vskip5pt

$(2)$ \ $\phi: P\otimes_A X_n\longrightarrow X'$ is injective;

\vskip5pt

$(3)$ \ $\cok \phi=(X_i/(\bigoplus\limits_{p\in \hp(n\rightarrow
i)}\Ima X_p), \ \widetilde{X_\alpha}, \ i\in Q'_0, \ \alpha\in
Q'_1)$, where for each $\alpha: j\longrightarrow i$ in $Q'_1$,
$$\widetilde{X_\alpha}: X_j/(\bigoplus\limits_{q\in \hp(n\rightarrow
j)}\Ima X_q) \longrightarrow {X_i}/ ({\bigoplus\limits_{p\in
\hp(n\rightarrow i)}\Ima X_p})$$ is the $A$-map induced by
$X_{\alpha}$.
\end{lem}

\noindent {\bf Proof.} By Lemma \ref{mono} and its proof it suffices
to prove $(1)$. For each $i\in Q'_0$, set $l_i = 0$ if
$\hp(n\rightarrow i)$ is empty, and $l_i = {\rm max} \{ \ l(p) \
\mid p\in \hp(n\rightarrow i)\}$ if otherwise, where $l(p)$ is the
length of $p$. We prove $(1)$ by using induction on $l_i$. If $l_i =
0$, then $(1)$ trivially holds. Suppose $l_i\geqslant 1$. Let
$\sum\limits_{p\in \hp(n\rightarrow i)} X_p(x_{n,p})=0$ for $x_{n,
p}\in X_n$. Since $${\sum\limits_{p \in \hp(n \to
i)- \ha(n \to i)}\Ima X_p} = \sum\limits_{\begin {smallmatrix} \alpha\in Q'_1\\
e(\alpha) = i \end{smallmatrix}} X_{\alpha} (\sum\limits_{q \in
\hp(n \to s({\alpha}))}\Ima X_q),$$ we have
\begin{align*} 0 = \sum\limits_{p\in \hp(n\rightarrow i)} X_p(x_{n,p})
& =\sum\limits_{\alpha \in  \ha(n\rightarrow i)} X_\alpha(x_{n,
\alpha}) +\sum\limits_{p \in \hp (n\rightarrow i)- \ha(n \to i)}
X_p(x_{n, p}) \\ & =\sum\limits_{\alpha \in  \ha(n\rightarrow i)}
X_\alpha(x_{n,\alpha})+ \sum\limits_{\begin {smallmatrix}\beta \in
Q'_1\\ e(\beta)=i
\end{smallmatrix}} X_\beta(\sum \limits_{q\in \hp(n\rightarrow
s(\beta))}X_q(x_{n,\beta q})),\end{align*} by $(m2)$ we know
$X_\alpha(x_{n,\alpha})=0$ for $ \alpha \in \ha(n\rightarrow i)$,
and $ X_\beta (\sum \limits_{q\in \hp(n\rightarrow s(\beta))}
X_q(x_{n, \beta q}))=0$ for $\beta \in Q'_1$ with $e(\beta)=i$. So
$\sum \limits_{q\in \hp(n\rightarrow s(\beta))} X_q(x_{n, \beta
q})=0$ by $(m1)$. Since $l_{s(\beta)} < l_i$ for each $\beta\in
Q_1'$ with $e(\beta) = i$, it follows from the inductive hypothesis
that $ X_q(x_{n, \beta q})=0$ for $\beta \in Q'_1, \ e(\beta)=i$ and
$q\in \hp(n\rightarrow s(\beta))$. This proves $(1)$. $\s$

\vskip10pt

\begin{lem} \label {cokerm2} \ Let $X=\left(\begin{smallmatrix}
X'\\
X_n
\end{smallmatrix}\right)_{\phi}$ be a monic $\Lambda$-module.
Then $\cok\phi$ is a monic $\Lambda'$-module.
\end{lem}

\vskip5pt

\noindent {\bf Proof.} \  We need to prove that for each $i\in
Q'_0$, the $\m'$-map
$$(\widetilde{X_\alpha})_{\alpha\in Q_1', \ e(\alpha)=i}: \ \bigoplus\limits_{\begin {smallmatrix}\alpha\in
Q'_1\\e(\alpha)=i\end{smallmatrix}}(X_{s(\alpha)}/(\bigoplus\limits_{q\in
\mathcal P (n\rightarrow s(\alpha))}{\rm Im}X_q))\longrightarrow
X_{i}/(\bigoplus\limits_{p\in \mathcal P (n\rightarrow i)}{\rm
Im}X_p)$$ is injective. For this, assume $\sum\limits_{\begin
{smallmatrix}\alpha\in Q'_1\\e(\alpha)=i\end{smallmatrix}}
 \widetilde{X_\alpha}(\overline{x_{{s(\alpha)}, \alpha}})=0,$
where $\overline {x_{s(\alpha),\alpha}}$ \ is the image of
$x_{s(\alpha), \alpha}\in X_{s(\alpha)}$ \ in
$X_{s(\alpha)}/(\bigoplus\limits_{q\in \mathcal P(n\rightarrow
s(\alpha))}{\rm Im}X_q).$ Then $\sum\limits_{\begin
{smallmatrix}\alpha\in Q'_1\\e(\alpha)=i\end{smallmatrix}}
X_\alpha(x_{{s(\alpha)},\alpha})\in \bigoplus\limits_{p\in
\hp(n\rightarrow i)}\Ima X_p$. So there are $x_{n,p}\in X_n$ such
that

$$\sum\limits_{\begin {smallmatrix}\alpha\in
Q'_1\\e(\alpha)=i\end{smallmatrix}} X_\alpha(x_{{s(\alpha)},\alpha})
=\sum\limits_{p\in \hp(n\rightarrow i)}X_p(x_{n,p}).$$ Thus
\begin{align*} 0& = \sum\limits_{\begin {smallmatrix}\alpha\in
Q'_1\\e(\alpha)=i\end{smallmatrix}}
X_{\alpha}(x_{{s(\alpha)},\alpha})- \sum\limits_{p\in
\hp(n\rightarrow i)}X_p(x_{n,p})\\ &=\sum\limits_{\begin
{smallmatrix}\alpha\in Q'_1\\e(\alpha)=i\end{smallmatrix}}
 X_{\alpha}(x_{{s(\alpha)},\alpha})
-\sum\limits_{\beta\in \ha(n\rightarrow i)}X_\beta (x_{n,\beta})
-\sum\limits_{\begin{smallmatrix}\alpha \in
Q'_1\\e(\alpha)=i\end{smallmatrix}}
X_{\alpha}(\sum\limits_{q\in \hp(n\rightarrow s(\alpha))}X_q(x_{n,\alpha q}))\\
&=\sum\limits_{\begin{smallmatrix}\alpha \in
Q'_1\\e(\alpha)=i\end{smallmatrix}}
X_\alpha(x_{s(\alpha),\alpha}-\sum\limits_{q\in \hp(n\rightarrow
s(\alpha))}X_q(x_{n,\alpha q}))- \sum\limits_{\beta \in
\ha(n\rightarrow i)}X_\beta(x_{n,\beta}).\end{align*}

\noindent Using the assumption on $X$ we get $x_{s(\alpha),\alpha}=
\sum\limits_{q \in \hp(n\rightarrow s(\alpha))}X_q(x_{n,\alpha q})$,
i.e., $\overline { x_{s(\alpha), \alpha}} =0$. $\s$

\vskip10pt

\begin{lem} \label{G4} \ Let $X=\left(\begin{smallmatrix}
X'\\
X_n
\end{smallmatrix}\right)_{\phi}$ be a monic $\Lambda$-module satisfying $(G)$.
Then $$(X_i/ (\bigoplus\limits_{p\in \hp(n\rightarrow i)}\Ima X_p))/
(\bigoplus \limits_{\begin{smallmatrix}\alpha \in
Q'_1\\e(\alpha)=i\end{smallmatrix}}\Ima \widetilde {X_\alpha})$$ is
a Gorenstein-projective $A$-module, $\forall \ i \in Q'_0$.
\end{lem}

\noindent{\bf Proof. } \  Since ${\bigoplus\limits_{p \in \hp(n \to
i)- \ha(n \to i)}\Ima X_p} \subseteq  \sum\limits_{\begin
{smallmatrix} \beta \in Q_1\\ e(\beta) = i \end{smallmatrix}}\Ima
X_{\beta}$, it follows that

\begin{align*} \ \ \ \ \ \ \ \ \ \ \ \  \ \ \ \ \ \ \ \ \ \ \ \
\sum\limits_{\begin {smallmatrix} \alpha\in Q'_1\\ e(\alpha) = i
\end{smallmatrix}}\Ima \widetilde{X_{\alpha}} & = ({\sum\limits_{\begin
{smallmatrix} \alpha\in Q'_1\\ e(\alpha) = i
\end{smallmatrix}}\Ima X_{\alpha} + \bigoplus\limits_{p \in \hp(n \to
i)}\Ima X_p})/ ({\bigoplus\limits_{p \in \hp(n \to i)}\Ima X_p})\\
& = ({\sum\limits_{\begin {smallmatrix} \beta \in Q_1\\
e(\beta) = i \end{smallmatrix}}\Ima X_{\beta} + \bigoplus\limits_{p
\in \hp(n \to i)-\ha(n \to i)}\Ima X_p})/ ({\bigoplus\limits_{p \in
\hp(n \to i)}\Ima X_p})\\ & = ({\sum\limits_{\begin {smallmatrix}
\beta \in Q_1\\ e(\beta) = i \end{smallmatrix}}\Ima
X_{\beta}})/({\bigoplus\limits_{p \in \hp(n \to i)}\Ima X_p}) \\ &
\stackrel{(m2)}=
({\bigoplus\limits_{\begin {smallmatrix} \beta \in Q_1\\
e(\beta) = i \end{smallmatrix}}\Ima
X_{\beta}})/({\bigoplus\limits_{p \in \hp(n \to i)}\Ima X_p}), \ \ \
\ \ \ \ \ \ \ \ \  \ \ \ \ \ \ \ \ \ \ \ \  \ \ \ \ \ \ \ \ \ \ \ \
(4.1)
\end{align*} and hence the desired quotient is
${X_i}/({\bigoplus\limits_{\begin {smallmatrix}\beta \in Q_1\\
e(\beta) = i \end{smallmatrix}}\Ima X_{\beta}})$, which is
Gorenstein-projective by $(G)$.   $\s$

\vskip10pt

\begin{lem}  \label{G3} \ Let $X=\left(\begin{smallmatrix}
X'\\
X_n
\end{smallmatrix}\right)_{\phi}$ be a monic $\Lambda$-module satisfying $(G)$.
Then for each $i\in Q_0'$, ${X_i}/ ({\bigoplus\limits_{p\in
\hp(n\rightarrow i)}\Ima X_p})$ is a Gorenstein-projective
$A$-module.
\end{lem}

\noindent {\bf Proof. } \ We prove the assertion by using induction
on $l_i$, which is defined in the proof of Lemma \ref{mono2}. If
$i\in Q_0'$ with $l_i = 0$, then the assertion follows from $(G)$.
\vskip5pt

Suppose $l_i\geqslant 1$. Since $\bigoplus \limits_{p\in
\hp(n\rightarrow i)}\Ima X_p
\subseteq\bigoplus\limits_{\begin{smallmatrix}\alpha\in
Q_1\\e(\alpha)=i\end{smallmatrix}}\Ima X_\alpha,$ we have the
following exact sequence
$$0\longrightarrow ({\bigoplus\limits_{\begin{smallmatrix}\alpha\in Q_1\\e(\alpha)=i\end{smallmatrix}}
\Ima X_\alpha})/({\bigoplus\limits_{p\in \hp(n\rightarrow i)}\Ima
X_p}) \longrightarrow {X_i}/({\bigoplus\limits_{p\in
\hp(n\rightarrow i)}\Ima X_p}) \longrightarrow
{X_i}/({\bigoplus\limits_{\begin{smallmatrix}\alpha\in
Q_1\\e(\alpha)=i\end{smallmatrix}}\Ima X_\alpha}) \longrightarrow
0,$$ by $(G)$ the term at the right hand side is
Gorenstein-projective. It suffices to prove that the term at the
left hand side is Gorenstein-projective. While by $(4.1)$ it is
$\bigoplus\limits_{\begin {smallmatrix} \alpha\in Q'_1\\ e(\alpha) =
i \end{smallmatrix}}\Ima \widetilde{X_{\alpha}}$. By Lemma \ref
{cokerm2} each $\widetilde{X_{\alpha}}$ is injective, it follows
that $\Ima \widetilde {X_\alpha}\cong {X_j}/({\bigoplus\limits_{p\in
\hp(n\rightarrow j)}\Ima X_p})$, where $j=s(\alpha)$. Since
$l_j<l_i$, it follows from the inductive hypothesis that
${X_j}/({\bigoplus\limits_{p\in \hp(n\rightarrow j)}\Ima X_p})$ is
Gorenstein-projective. This completes the proof.  $\s$

\vskip10pt

\begin{lem} \label{sufficiency} The sufficiency in Theorem \ref{mainthm} holds. That is, if
$X = (X_i, \ X_\alpha, \ i\in Q_0, \ \alpha\in Q_1)$ is a monic
$\m$-module satisfying $(G)$, then $X$ is Gorenstein-projective.
\end{lem}

\noindent {\bf Proof.} \ Using induction on $n=| Q_0|$. The
assertion clearly holds for $n=1$. Suppose that the assertion holds
for $n-1$ with $n\geq 2$. It suffices to prove that $X$ satisfies
the conditions $(i)$, $(ii)$ and $(iii)$ in Lemma \ref{mainlem}.

\vskip5pt

The condition $(i)$ is contained in $(G)$; and the condition $(ii)$
follows from Lemma \ref{mono2}$(2)$. By Lemma \ref{cokerm2}
$\cok\phi$ is a monic $\Lambda'$-module; and by Lemmas \ref{G4} and
\ref{G3} we know that $\cok\phi$ satisfies $(G)$. It follows from
the inductive hypothesis that the condition $(iii)$ is satisfied.
$\s$

\vskip10pt

\begin{lem} \label{X'gp} Let $X = (X_i, \ X_\alpha, \ i\in Q_0, \ \alpha\in Q_1)$ be a
$\Lambda$-module with $X_n$ a Gorenstein-projective $A$-module. Then
$P\otimes_A X_n$ is a Gorenstein-projective $\Lambda'$-module, where
$P$ is defined in 4.2.
\end{lem}

\noindent {\bf Proof.} \ Let $P(n)$ be the indecomposable projective
$kQ$-module at vertex $n$. Writing ${\rm rad} P(n)$ as a
representation of $Q'$ over $k$, we have ${\rm rad} P(n) = (k^{m_i},
\ f_\alpha, \ i\in Q'_0, \ \alpha\in Q'_1),$ where
$m_i=|\mathcal{P}(n \to i)|$ for each $i \in Q'_0$. By the
construction of $P(n)$ we know that ${\rm rad} P(n)$ has the
following three properties:

\vskip5pt

$(1)$ \ each $f_\alpha: k^{m_{s(\alpha)}}\longrightarrow
k^{m_{e(\alpha)}}$ is injective;

\vskip5pt

$(2)$ \ for each $i\in Q_0'$  there holds $\sum\limits_{\begin
{smallmatrix} \alpha\in Q'_1\\ e(\alpha) = i
\end{smallmatrix}}\Ima
f_\alpha = \bigoplus\limits_{\begin {smallmatrix} \alpha\in Q_1'\\
e(\alpha) = i \end{smallmatrix}}\Ima f_\alpha;$

\vskip5pt

$(3)$ \ for each $i \in Q'_0$,  $k^{m_i}/ ({\bigoplus\limits_{\begin {smallmatrix} \alpha\in Q'_1\\
e(\alpha) = i \end{smallmatrix}}\Ima f_\alpha}) \cong
k^{|\mathcal{A}(n \to i)|}$ as $k$-spaces.

\vskip5pt

It follows that
$$P\otimes_A X_n \cong  ({\rm rad} P(n)\otimes_kA)\otimes_A X_n \cong {\rm rad} P(n)\otimes_kX_n = (X_n^{m_i}, \ f_\alpha\otimes_k {\rm id}_{X_n}, \ i\in Q'_0, \
\alpha\in Q'_1).$$ By $(1), (2)$ and $(3)$ we clearly see that
$P\otimes_A X_n$ is a monic $\Lambda'$-module satisfying $(G)$ (for
example, by $(3)$ we know that
$X_n^{m_i}/({\bigoplus\limits_{\begin {smallmatrix} \alpha\in Q'_1\\
e(\alpha) = i \end{smallmatrix}}\Ima (f_\alpha}\otimes_k{\rm
id}_{X_n})) \cong X_n^{|\mathcal{A}(n \to i)|}$ is a
Gorenstein-projective $A$-module).  Now the assertion follows from
Lemma \ref{sufficiency}. $\s$

\vskip10pt

\subsection{} \ {\bf Proof of Theorem 4.1} \ By Lemma \ref{sufficiency} it
remains to prove the necessity, i.e., if $X$ is a
Gorenstein-projective $\Lambda$-module, then $X$ is a monic
$\Lambda$-module satisfying $(G)$. Using induction on $n=|Q_0|$. The
assertion is clear for $n=1$. Suppose that the assertion holds for
$n-1$ with $n\geq 2$. We write as $X=\left(\begin{smallmatrix}
X'\\
X_n
\end{smallmatrix}\right)_{\phi}$. Then $X$ satisfies the conditions $(i), (ii)$ and
$(iii)$ in Lemma \ref{mainlem}.

\vskip5pt

By the condition $(i)$ and Lemma \ref{X'gp} we know that $P\otimes_A
X_n$ a Gorenstein-projective $\Lambda'$-module. Then by the
conditions $(ii)$ and $(iii)$ we know that $X'\in\mathcal{GP}(\m')$
since $\mathcal{GP}(\m')$ is closed under extensions. By the
inductive hypothesis $X'$ is a monic $\Lambda'$-module satisfying
$(G)$, thus the following properties hold:

\vskip5pt

$(1)$ \ $X_{\beta}$ is injective for each $\beta \in Q'_1$; and

\vskip5pt

$(2)$ \ $X_{i}$ is Gorenstein-projective for each $i \in Q'_0$.

\vskip5pt

By $(1)$, the condition $(ii)$ and Lemma \ref{mono} we know that

\vskip5pt

$(3)$ \ $X_\alpha$ is injective for each $\alpha \in Q_1$; and

\vskip5pt

$(4)$ \ $\sum\limits_{p \in \hp(n \to i)}\Ima X_p
=\bigoplus\limits_{p \in \hp(n \to i)}\Ima X_p, \ \forall \ i\in
Q'_0. $

\vskip5pt

Since $\cok \phi=({X_i}/({\bigoplus\limits_{p \in \hp(n \to i)}\Ima
X_p}), \ \widetilde {X_{\alpha}}, \ i\in Q'_0, \ \alpha\in Q_1')$ is
a Gorenstein-projective $\Lambda'$-module, it follows from the
inductive hypothesis that the following properties hold:

\vskip5pt

$(5)$ \  for each $\alpha\in Q'_1$, $\widetilde{X_{\alpha}}$ is
injective; and

\vskip5pt

$(6)$ \ $\sum\limits_{\begin {smallmatrix} \alpha\in Q'_1\\
e(\alpha) = i \end{smallmatrix}}\Ima \widetilde{X_\alpha} =
\bigoplus\limits_{\begin {smallmatrix} \alpha\in Q'_1\\ e(\alpha) =
i
\end{smallmatrix}}\Ima \widetilde{X_\alpha}, \ \forall \ i\in Q_0'.$

\vskip10pt

We first prove Claim 1: $X$ satisfies $(m2)$. In fact, suppose
$$\sum\limits_{\begin {smallmatrix} \alpha\in Q_1\\ e(\alpha) = i
\end{smallmatrix}}X_{\alpha}(x_{s{(\alpha)},\alpha}) =0. \eqno (*)$$ Since
$$\sum\limits_{\begin {smallmatrix}
\alpha\in Q_1\\ e(\alpha) = i
\end{smallmatrix}}X_{\alpha}(x_{s{(\alpha)},\alpha}) =
\sum\limits_{\begin {smallmatrix} \alpha\in \ha(n\to i)
\end{smallmatrix}}X_{\alpha}(x_{s{(\alpha)}, \alpha}) +
\sum\limits_{\begin {smallmatrix} \alpha\in Q'_1\\ e(\alpha) = i
\end{smallmatrix}}X_{\alpha}(x_{s{(\alpha)},\alpha}),$$ it follows
that

\begin{align*} \sum\limits_{\begin {smallmatrix} \alpha\in Q'_1\\
e(\alpha) = i
\end{smallmatrix}}\widetilde{X_{\alpha}}(\overline{x_{s{(\alpha)},\alpha}})
& = \sum\limits_{\begin {smallmatrix} \alpha\in Q'_1\\ e(\alpha) = i
\end{smallmatrix}}X_{\alpha}(x_{s{(\alpha)},\alpha})+
(\bigoplus\limits_{p\in \hp(n\to i)}\Ima X_p)
\\ & \stackrel {\mbox {by} (*)} = -\sum\limits_{\begin {smallmatrix} \alpha\in \ha(n\to i)
\end{smallmatrix}}X_{\alpha}(x_{s{(\alpha)},\alpha}) +
(\bigoplus\limits_{p\in \hp(n\to i)}\Ima X_p) =0.\end{align*} Then
by $(6)$ we have $\widetilde{X_{\alpha}}(\overline{x_{s(\alpha),
\alpha}})=0$; and by $(5)$ we know $\overline{x_{s(\alpha),
\alpha}}=0$ for each $\alpha \in Q'_1$ with $e(\alpha)=i$. This
means that there are $x_{n, q}\in X_n$ such that
$$x_{s(\alpha),\alpha}=
\sum\limits_{q\in \hp(n \to s(\alpha))}X_{q}(x_{n,q})\in \sum\limits_{q\in \hp(n \to s(\alpha))}\Ima X_q$$
for each  $\alpha \in Q'_1$ with $e(\alpha)=i$. By $(*)$ we have
$$0=\sum\limits_{{\alpha}\in \ha(n \to i)}X_{\alpha}(x_{n,{\alpha}})
+\sum\limits_{\begin {smallmatrix} \alpha\in Q'_1\\ e(\alpha) = i
\end{smallmatrix}}X_{\alpha}( \sum\limits_{q\in \hp(n \to
s(\alpha))}X_{q}(x_{n,q})).$$ By $(4)$ we know that
$X_{\alpha}(x_{n,{\alpha}})=0, \forall \alpha \in \ha(n \to i)$, and
$X_{\alpha}X_q(x_{n, q})=0, \forall \alpha\in Q'_1$ with $e(\alpha)
= i$ and $q\in \hp(n \to s(\alpha))$. Thus
$X_{\alpha}(x_{s{(\alpha)},\alpha})=0, \forall \alpha\in Q_1$ with
$e(\alpha) = i$. This proves Claim 1.

\vskip5pt

We now prove Claim 2:
$X_i/(\bigoplus\limits_{\begin {smallmatrix}\beta \in Q_1\\
e(\beta) = i \end{smallmatrix}}
 \Ima X_\beta)$ is a Gorenstein-projective $A$-module for each $i\in Q_0$.

\vskip5pt\noindent In fact, since $\cok \phi$ is a
Gorenstein-projective $\m'$-module, by the inductive hypothesis we
know that $$(X_i/(\bigoplus\limits_{p \in \hp(n
\to i)}\Ima X_p))/(\bigoplus\limits_{\begin {smallmatrix} \alpha\in Q'_1\\
e(\alpha) = i \end{smallmatrix}} \Ima \widetilde {X_\alpha})$$ is a
Gorenstein-projective $A$-module: it is exactly the desired module
by $(4.1)$.

\vskip10pt

Now, $(3)$ and Claim 1 mean that $X$ is a monic $\Lambda$-module;
and $(2)$, the condition $(i)$, together with Claim 2 mean that $X$
satisfies $(G)$. This completes the proof. $\s$

\section {\bf Applications}

\vskip10pt

We include some applications of Theorem \ref{mainthm}. Let $\m$ be
the path algebra of finite acyclic quiver $Q$ over
finite-dimensional algebra $A$. Recall that ${\rm Mon}(Q, A)$
denotes the full subcategory of ${\rm Rep}(Q, A)$ consisting of the
monic representations of $Q$ over $A$.

\vskip10pt

\subsection{} As a consequence of Theorem \ref{mainthm}, we get the following characterization of
self-injectivity.

\vskip10pt

\begin{thm} \label{self} \ Let $A$ be a finite-dimensional algebra. Then the following are equivalent:

\vskip5pt

$(i)$ \  $A$ is self-injective;

\vskip5pt

$(ii)$ \ For any  finite acyclic quiver $Q$, there holds $\mathcal
{GP}(A\otimes_kkQ) = {\rm Mon}(Q, A)$;

\vskip5pt

$(iii)$ \ There is a finite acyclic quiver $Q$, such that $\mathcal
{GP}(A\otimes_kkQ) = {\rm Mon}(Q, A)$.
\end{thm}

\noindent{\bf Proof.} \ $(i) \Longrightarrow (ii)$: If $A$ is
self-injective, then every $A$-module is Gorenstein-projective, and
hence $(ii)$ follows from Theorem \ref{mainthm}. The implication
$(ii)\Longrightarrow (iii)$ is clear.

$(iii) \Longrightarrow (i)$: Take a sink of $Q$, say vertex $1$, and
consider the representation $X$ of $Q$ over $A$, where $X_1 =
\Hom_A(A, k)$ and $X_i = 0$ if $i\ne 1$. Then $X$ is a monic
$\m$-module, and hence by assumption it is Gorenstein-projective. So
we have a complete $\m$-projective resolution
$$\cdots \longrightarrow P^{-1}\longrightarrow P^{0}
\stackrel{d^0}{\longrightarrow} P^{1}\longrightarrow \cdots $$ with
$X\cong \operatorname{Ker}d^0$. By taking the $1$-st branch we get
an exact sequence $$\cdots \longrightarrow P^{-1}_1\longrightarrow
P^{0}_1 \stackrel{d^0_1}{\longrightarrow} P^{1}_1\longrightarrow
\cdots $$ with $\Ker d^0_1\cong \Hom_A(A, k)$. Note that each
$P^i_1$ is a projective $A$-module. Thus injective $A$-module
$\Hom_A(A, k)$ is projective, i.e., $A$ is self-injective.  $\s$

\vskip10pt

\begin{exm}  Taking $Q = \underset n \bullet \longrightarrow \cdots \longrightarrow
 \underset 1\bullet$ in Theorem \ref{self} we get: $A$ is a self-injective algebra if and only if the
Gorenstein-projective $T_n(A)$-modules are exactly the monic
$T_n(A)$-modules. Under the assumption that $A$ is Gorenstein, this
result was obtained in Theorem 4.4 of [Z].
\end{exm}

\vskip10pt

Let $D^b(\m)$ be the bounded derived category of $\m$, and
$K^b(\mathcal P(\m))$ the bounded homotopy category of $\mathcal
P(\m)$. By definition the singularity category $D^b_{sg}(\m)$ of
$\m$ is the Verdier quotient $D^b(\m)/K^b(\mathcal P(\m))$. If $\m$
is Gorenstein, then there is a triangle-equivalence
$D^b_{sg}(\m)\cong
 \underline {\mathcal
{GP}(\m)},$ where $ \underline {\mathcal {GP}(\m)}$ is the stable
category of \ ${\mathcal {GP}(\m)}$ modulo $\mathcal P(\m)$ (see
[Hap], Theorem 4.6; also [Buch], Theorem 4.4.1). Note that if $A$ is
Gorenstein, then $\m = A\otimes_kkQ$ is Gorenstein, by Proposition
2.2 in [AR], which claims that $A\otimes_kB$ is Gorenstein if and
only if $A$ and $B$ are Gorenstein. So we have

\vskip10pt

\begin{cor} \label{cor} \ Let $A$ be a finite-dimensional
Gorenstein algebra, and $Q$ a finite acyclic quiver. Let $\m =
A\otimes_kkQ$. Then there is a triangle-equivalence \
$D^b_{sg}(\m)\cong  \underline {\mathcal {GP}(\m)}.$
\end{cor}

\vskip10pt

\subsection{} Before giving the next application we recall the tensor product of two finite
quivers. Let $Q$ and $Q'$ be finite quivers (not necessarily
acyclic). By definition the tensor product $Q\otimes Q'$ is the
quiver with
$$(Q\otimes Q')_0 = Q_0\times Q_0', \ \  \mbox{and} \ \ (Q\otimes Q')_1 = (Q_1\times
Q_0') \bigcup (Q_0\times Q_1').$$ More explicitly, if $\alpha:
i\longrightarrow j$ is an arrow of $Q$, then for each vertex $t'\in
Q_0'$ there is an arrow $(\alpha, t'): (i, t') \longrightarrow (j,
t')$ of $Q\otimes Q'$; and if $\beta': s'\longrightarrow t'$ is an
arrow of $Q'$, then for each vertex $i\in Q_0$ there is an arrow
$(i, \beta'): (i, s') \longrightarrow (i, t')$ of $Q\otimes Q'$.

\vskip10pt

Let $A = kQ/I$ and $B = kQ'/I'$ be two finite-dimensional
$k$-algebra, where $Q$ and $Q'$ are finite quivers (not necessarily
acyclic), $I$ and $I'$ are  admissible ideals of $kQ$ and $kQ'$,
respectively. Then
$$A\otimes_k B \cong k(Q\otimes Q')/I\Box I',$$
where $I\Box I'$ is the ideal of $k(Q\otimes Q')$ generated by
$(I\times Q_0')\bigcup (Q_0\times I')$ and the following elements
$$(\alpha, t')(i, \beta')-(j, \beta')(\alpha, s'),$$
where $\alpha: i\longrightarrow j$ is an arrow of $Q$, and $\beta':
s'\longrightarrow t'$ is an arrow of $Q'$. See for example [L]. Note
that $I\Box I'$ may not be zero even if $I = 0 = I'$. Thus we have
the following

\vskip10pt

\noindent {\bf Fact:} \   $A\otimes_kB$ is hereditary (i.e., $I\Box
I' = 0$) if and only if either

\vskip5pt

$(i)$ \ $A \cong  k^{|Q_0|}$ as algebras, and $I' = 0$;  or

\vskip5pt

$(ii)$ \ $B \cong k^{|Q_0'|}$ as algebras, and $I = 0$.

\vskip10pt

\subsection{} We describe when $\m$ is hereditary via monic $\m$-modules.

\vskip10pt

\begin{thm} \label{hereditary} \ Let $\m$ be the path algebra of finite quiver $Q$ over
$A$, where $Q$ is  acyclic with $|Q_1| \ne 0$, and $A$ is a
finite-dimensional basic algebra over an algebraically closed field
$k$. Then $\mathcal P(\m) = {\rm Mon}(Q, A)$ if and only if $\m$ is
hereditary.
\end{thm}

\noindent{\bf Proof.} Without loss of generality we may assume that
$A$ is connected (an algebra is connected if it can not be a product
of two non-zero algebras).

If $\m = A\otimes_kkQ$ is hereditary, then by the fact above and the
assumption of $Q$ we have $A = k$, and hence ${\rm Mon}(Q, k) =
\mathcal {GP}(kQ)$ by Theorem \ref{mainthm}. It follows that
$${\rm Mon}(Q, A) = {\rm Mon}(Q, k) = \mathcal {GP}(kQ) = \mathcal
P(kQ) = \mathcal P(\m).$$

Conversely, if $A\ne k$, then $A$ is not semi-simple since $A$ is
assumed to be connected and basic and $k$ is assumed to be
algebraically closed. It follows that there is a non-projective
$A$-module $M$. Take a sink of $Q$, say vertex $1$, and consider
$\m$-module $X = M\otimes_k P(1)$, where $P(1)$ is the simple
projective $kQ$-module at vertex $1$. Then as a representation of
$Q$ over $A$ we have $X = (X_i, \ i\in  Q_0)$ with $X_1 = M$ and
$X_i = 0$ for $i\ne 1$. It is clear that $X\in {\rm Mon}(Q, A)$, but
$X\notin \mathcal P(\m).$ $\s$


\vskip15pt



\begin{thebibliography}{99}

\bibitem[AB]{AB} M. Auslander, M. Bridger, Stable module
theory, Mem. Amer. Math. Soc. 94., Amer. Math. Soc., Providence,
R.I., 1969.

\bibitem[AM] {AM} L. L. Avramov, A. Martsinkovsky, Absolute, relative, and Tate cohomology of modules of
finite Gorenstein dimension, Proc. London Math. Soc. 85(3)(2002),
393-440.

\bibitem[AR]{AR} M. Auslander and I. Reiten, Cohen-Macaulay and Gorenstein artin algebras, In: Representation
theory of finite groups and finite-dimensional algebras (Proc. Conf.
at Bielefeld, 1991), 221-245, Progress in Math. vol. 95,
Birkh\"auser, Basel, 1991.

\bibitem[ARS]{ARS} M. Auslander, I. Reiten, S. O.
Smal${\o}$, Representation Theory of Artin Algebras, Cambridge
Studies in Adv. Math. 36., Cambridge Univ. Press, 1995.

\bibitem[B]{B}
A. Beligiannis, Cohen-Macaulay modules, (co)torsion pairs and
virtually Gorenstein algebras, J. Algebra 288(1)(2005), 137-211.

\bibitem[Buch]{B} R.-O. Buchweitz, Maximal Cohen-Macaulay modules and Tate cohomology over Gorenstein rings,
Unpublished manuscript, Hamburg (1987), 155pp.

\bibitem[C]{C} X. W. Chen, The stable monomorphism category
of a Frobenius category, Math. Res. Lett. 18(1)(2011), 125-137.

\bibitem[EJ]{EJ} E. E. Enochs, O. M. G. Jenda, Relative homological
algebra, De Gruyter Exp. Math. 30.
Walter De Gruyter Co., 2000.

\bibitem[GZ]{GZ} N.
Gao, P. Zhang, Gorenstein derived categories, J. Algebra 323(2010),
2041-2057.

\bibitem[Hap]{Hap} D. Happel, On Gorenstein algebras, in:
Representation theory of
finite groups and finite-dimensional algebras, Prog. Math. 95,
389-404, Birkh\"user, Basel, 1991.

\bibitem[Hol]{Hol} H. Holm, Gorenstein
homological dimensions, J. Pure Appl. Algebra 189(1-3)(2004),
167-193.

\bibitem[IKM]{IKM} O.
Iyama,  K. Kato, J. I. Miyachi, Recollement on homotopy categories
and Cohen-Macaulay modules, available in arXiv: math. RA
0911.0172.

\bibitem[KLM]{KLM} D. Kussin, H. Lenzing, H. Meltzer,
Nilpotent operators and weighted projective lines, avaible in arXiv:
math. RT 1002.3797.

\bibitem[L]{L} Z. Leszczy\'nski, On the representation type of tensor product of algebras, Fundamenta Math. 144(1994),
143-161.

\bibitem[LZ]{LZ} Z. W. Li, P. Zhang, A construction of
Gorenstein-projective modules, J. Algebra 323(2010),  1802-1812.



\bibitem [R]{R} C. M. Ringel, Tame algebras and integral
quadratic forms, Lecture Notes in Math. 1099, Springer-Verlag, 1984.

\bibitem[RS1]{RS1} C. M. Ringel, M. Schmidmeier, Submodules
categories of wild representation type, J. Pure Appl. Algebra
205(2)(2006), 412-422.
\bibitem[RS2]{RS2}
C. M. Ringel, M. Schmidmeier, The Auslander-Reiten translation in
submodule categories, Trans. Amer. Math. Soc. 360(2)(2008), 691-716.
\bibitem[S]{S} D. Simson, Representation types of the category of subprojective
representations of a finite poset over $K[t]/(t^m)$ and a solution
of a Birkhoff type problem, J. Algebra 311(2007),
1-30.

\bibitem[XZ]{XZ} B. L. Xiong, P. Zhang, Gorenstein-projective modules
over triangular matrix Artin algebras, to appear in: J. Algebra
Appl.
\bibitem[Z]{Z} P. Zhang, Monomorphism categories, cotilting
theory, and Gorenstein-projective modules, J. Algebra 339(2011),
180-202.
\end{thebibliography}
\end{document}